\documentclass[twoside,
]{amsart}

  \usepackage[latin1]{inputenc}
  \usepackage{amsmath}
  \usepackage{amssymb}
  \usepackage{amsthm}
  \usepackage{mathrsfs}
  \usepackage{bbm}
  \usepackage{pifont}
  \usepackage{graphicx}
  \usepackage{color}

	\newtheoremstyle{slanted}
	{}
	{}
	{\slshape}
	{}
	{\bfseries}
	{.}
	{ }
	{}
	
	\theoremstyle{slanted}
	\newtheorem{theo}{Theorem}[section]

	\newtheorem{question}[theo]{Question}
	\newtheorem{lemma}[theo]{Lemma}
	\newtheorem{definition}[theo]{Definition}
	\newtheorem{corollary}[theo]{Corollary}
	\newtheorem{example}[theo]{Example}

	\def\ind#1{\mathbbmss{1}_{#1}}
	\def\egdef{:=}

	\newcommand{\tend}[2]{\xrightarrow[#1\to#2]{}}

	\def\ind#1{\mathbbmss{1}_{#1}}
	\newcommand{\ZZ}{\mathbb{Z}}
	\renewcommand{\AA}{\mathbb{A}}
	\newcommand{\QQ}{\mathbb{Q}}
	\newcommand{\RR}{\mathbb{R}}
	\newcommand{\PP}{\mathbb{P}}
	\newcommand{\NN}{\mathbb{N}}
	\newcommand{\C}{\mathscr{C}}

\title[Pairwise-independent Joinings]{A Class of pairwise-independent Joinings}

\author{Élise Janvresse and Thierry de la Rue}

\dedicatory{Dedicated to the memory of Professor José de Sam Lazaro}

\address{Laboratoire de Math\'ematiques Rapha\"el Salem\\
	UMR 6085 CNRS -- Universit\'e de Rouen\\
	Avenue de l'Universit\'e\\
	B.P. 12\\
	F76801 Saint-\'Etienne-du-Rouvray Cedex}
\email{Elise.Janvresse@univ-rouen.fr, Thierry.de-la-Rue@univ-rouen.fr}

\begin{document}
\bibliographystyle{amsplain}

\keywords{pairwise-independent joining, subshift of finite type, multifold mixing}
\subjclass[2000]{37A35, 37B10}

\begin{abstract}
We introduce a special class of pairwise-independent self-joinings for a stationary process: Those for which one coordinate is a continuous function of the two others. We investigate which properties on the process the existence of such a joining entails. In particular, we prove that if the process is aperiodic, then it has positive entropy. Our other results suggest
that such pairwise independent, non-independent self-joinings exist only in very specific situations: Essentially when the process is a subshift of finite type topologically conjugate to a full-shift. This provides an argument in favor of the conjecture that 2-fold mixing implies 3-fold-mixing.
\end{abstract}

\maketitle

\section{Introduction}
The purpose of this article is to study in which cases it is possible to construct some special self-joinings for a stationary process.
Throughout the text, $\xi=(\xi_i)_{i\in\ZZ}$ is a stationary process taking values in the finite alphabet $\AA$, and $\mu$ is its probability distribution on $\AA^\ZZ$. For any $i<j$, we denote by $\xi_i^j$ the sequence of random variables $(\xi_i,\xi_{i+1},\ldots,\xi_j)$.

\begin{definition}
Let $k\ge2$. A \emph{$k$-fold self-joining} of $\xi$ is a stationary probability distribution on $(\AA^k)^\ZZ$ (canonically identified with $(\AA^\ZZ)^k$), whose $k$ marginals on $\AA^\ZZ$ are all equal to $\mu$.
\end{definition}

In the sequel, we will mostly concentrate on 2-fold and 3-fold self-joinings. A 3-fold self-joining of $\xi$ can be viewed as the joint distribution of three processes $\xi$, $\xi'$ and $\xi''$ with the same distribution $\mu$, provided this joint distribution is stationary.

\subsection{Pairwise-independent joinings}
\label{Sec:intro}
\begin{definition}
We call \emph{pairwise-independent joining of $\xi$} a 3-fold self-joining of $\xi$ such that the three 2-dimensional marginals are equal to $\mu\otimes\mu$.
\end{definition}

An obvious example of a pairwise-independent joining of $\xi$ is the product measure $\mu\otimes\mu\otimes\mu$, corresponding to the case where $\xi$, $\xi'$ and $\xi''$ are independent.

Two other classical examples arise for particular processes $\xi$. 
\begin{example}
\label{Bernoulli}
When $\AA=\ZZ/n\ZZ$ and $\xi$ is a uniform Bernoulli process (\textit{i.e.} the random variables $\xi_i$ are i.i.d. and uniformly distributed in $\AA$), we can take $\xi'$ as an independent copy of $\xi$, and set for all $i\in\ZZ$
$$ 
\xi''_i\egdef\xi_i+\xi'_i\quad\mod n.
$$
Then $\xi''$ follows the same distribution as $\xi$ and $\xi'$, is independent from $\xi$ and from $\xi'$, but the joint distribution is clearly not the product measure, since $\xi''$ is determined by $\xi$ and $\xi'$.
\end{example}
\begin{example}
\label{periodic}
When $\AA=\ZZ/n\ZZ$ for an odd integer $n$, and $\xi$ is a periodic process satisfying $\xi_{i+1}=\xi_i+1\mod n$ for all $i$, we can also take an independent copy $\xi'$ of $\xi$, and set
for all $i\in\ZZ$
$$ 
\xi''_i\egdef 2\xi_i-\xi'_i\quad\mod n.
$$
The same conclusion follows for the joint distribution of $\xi$, $\xi'$ and $\xi''$.
\end{example}

Whether there exist examples of a different nature is an oustanding problem in ergodic theory, which can be stated as follows.

\begin{question}
\label{question}
Does there exist a \emph{zero-entropy}, \emph{weakly mixing}, stationary process $\xi$ with a pairwise-independent joining $\lambda$,  $\lambda\neq\mu\otimes\mu\otimes\mu$?
\end{question}

Observe that the two preceding examples do not belong to the category of zero-entropy, weakly mixing processes. The above  open question is related to two important problems in ergodic theory: Whether 2-fold mixing implies 3-fold mixing, and whether 2-fold minimal self-joinings implies 3-fold minimal self-joinings. Indeed, a negative answer to Question~\ref{question} would prove both implications, since any counterexample to one of these problems would automatically lead to a non-trivial pairwise independent joining.

Many important results have already been obtained on this question. Most of them consider some special category of stationary processes (\textit{e.g.} finite-rank processes~\cite{kalikow1984, ryzhikov1993}, or processes with singular spectrum~\cite{host1991}), and prove that in this category a pairwise-independent 3-fold self-joining has to be the product measure. 
In the present work, we have a slightly different approach: Consider a stationary process which admits some kind of pairwise-independent 3-fold self-joining which is not the product measure, and see which other properties on the process this assumption entails.

\subsection{The PIJ property}

The special kind of pairwise-independent joining we want to study is the case where the third coordinate $\xi''$ is determined in a continuous way by the first two.

\begin{definition}
We say that the process $\xi$ has the \emph{PIJ property} if there exist a pairwise-independent joining $\lambda$ of $\xi$ and a continuous function $\varphi:\AA^\ZZ\times\AA^\ZZ\to \AA$ satisfying
\begin{equation}
\label{defPIJeq}
\xi''_0 = \varphi(\xi,\xi')\qquad(\lambda\mbox{-a.e.})
\end{equation}
\end{definition}

This amounts to saying that there exists some $p\in\NN$ such that 
\begin{equation}
\label{width}
\xi''_0 = \varphi(\xi_{-p}^{p},{\xi'}_{-p}^{p})\qquad(\lambda\mbox{-a.e.})
\end{equation}
We will sometimes refer to this integer $p$ as the \emph{width} of the joining $\lambda$.

The two examples given above (uniform Bernoulli and periodic process) obviously have the PIJ property, with $p=0$. It is straightforward to check that the PIJ property is preserved by topological conjugacy and independent product. Therefore, we get from these examples a whole family of PIJ processes, which all are subshifts of finite type. This leads to the following question.

\begin{question}
\label{question2}
Do there exist other processes enjoying the PIJ property than the processes arising from Examples~\ref{Bernoulli} and~\ref{periodic}?
\end{question}

A connected problem was raised by Ryzhikov \cite{ryzhikov1992} who asked under which condition it is possible to find a process $\zeta = \varphi(\xi,\xi')$ independent from $\xi$ and from $\xi'$ (where $\xi$ and $\xi'$ are two independent copies of the same process).

\subsection{Results}
Our first result, presented in Section~\ref{general}, states that a process which has the PIJ property is either periodic or has positive entropy. 
Therefore, no process possessing the PIJ property can answer positively to Question~\ref{question}. This generalizes the results obtained in~\cite{delarue2006b,delarue2006} in the case where the width of the joining is 0. Note that the first result in this direction was obtained by Ryzhikov in 1993: de Sam Lazaro, Thouvenot and Ryzhikov noticed that Assertion~4.4 in~\cite{ryzhikov1993a} could be reinterpreted in the following elegant form:
If $\xi$, $\xi'$ and $ \xi+\xi'$  are pairwise independent 2-valued processes, then $\xi$ is Bernoulli or periodic.

In Section~\ref{SSFT}, we study irreducible aperiodic subshifts of finite type having the PIJ property. 
We show that their Perron value is an integer, which implies that they are measure-theoretically isomorphic to a uniform Bernoulli shift (cf. Example~\ref{Bernoulli}). We even prove a stronger result: any irreducible aperiodic subshift of finite type with the PIJ property is shift-equivalent
\footnote{Williams \cite{williams1973} conjectured that shift equivalence was the same as topological conjugacy, which was disproved by Kim and Roush \cite{kim1997}. Even so, Williams' conjecture could be true when one of the subshift is the full shift.} to a full shift. 

Finally, we ask whether any process with the PIJ property is a subshift of finite type. 
We prove in Section~\ref{PIJimplySFT?} that a slightly stronger property than PIJ is sufficient: If, in the definition of the PIJ property, we further require that the first coordinate be a continuous function of $\xi'$ and $\xi''$, then the process is a subshift of finite type. 

All these results lead to the conjecture that there are no other processes with the PIJ property than those which arise from our two examples by topological conjugacy or direct product.

\section{General consequences of the PIJ property}
\label{general}
\subsection{Zero-entropy implies periodicity}

\begin{theo}
\label{positive entropy}
If $\xi$ has the PIJ property, then
\begin{itemize} 
\item either $\xi$ has positive entropy,
\item or $\xi$ is periodic.
\end{itemize}
\end{theo}

One of the main ingredients for the proof of this theorem is the property of \emph{quasi-uniformity} of the measure $\mu$ which is a consequence of the PIJ property. 

We call \emph{cylinder of length $\ell$} any event of the form $(\xi_i=a_i,\ldots,\xi_{i+\ell-1}=a_{i+\ell-1})$, provided that the sequence $(a_i,\ldots,a_{i+\ell-1})$ is seen with positive probability. This cylinder will be denoted by $[a_i^{i+\ell-1}]$.
The set of all words of length $\ell$ appearing with positive probability will be denoted by $\mathcal{L}_\ell$.

\begin{definition} The probability distribution $\mu$ is said to be \emph{quasi-uniform} if there exists a constant $K>0$ such that, if $B$ and $C$ are two cylinders of the same length,
\begin{equation}
\label{def_quasi_uniform}\dfrac{1}{K} \mu(B) \le \mu(C) \le K\mu(B).
\end{equation}

\end{definition}

\begin{lemma}
\label{lemma_unif}
Assume that $\xi$ has the PIJ property with some joining $\lambda$ of width $p$. 
\begin{itemize}
\item[(i)] If $A$ and $B$ are two cylinders with respective lengths $\ell\ge1$ and $\ell+2p$, then $\mu(A)\ge\mu(B)$.
\item[(ii)] $\mu$ is quasi-uniform.
\item[(iii)] $\mu$ has maximal entropy among all invariant probability measures with the same support.
\end{itemize}
\end{lemma}

\begin{proof}
Since $\mu$ is shift-invariant, there is no loss of generality in assuming that $A$ is of the form $[a_0^{\ell-1}]$ and $B=[b_{-p}^{p+\ell-1}]$. Under $\lambda$, $\xi'$ and $\xi''$ are independent, therefore 
$$
\lambda(\xi'\in B,\xi''\in A)>0.
$$
We can thus find some $\theta>0$ small enough such that
$$
C_\theta\egdef\{\xi\in\AA^\ZZ:\ \lambda(\xi'\in B,\xi''\in A|\xi)>\theta\}
$$
has positive probability. By \eqref{width} and by definition of $C_\theta$, 
$$
(\xi'\in B, \xi\in C_\theta)\subset (\xi''\in A, \xi\in C_\theta).
$$
Using pairwise independence, we get
$$
\mu(B)\mu(C_\theta) = \lambda(\xi'\in B, \xi\in C_\theta) \le \lambda (\xi''\in A, \xi\in C_\theta) = \mu(A)\mu(C_\theta).
$$
Hence $\mu(B)\le\mu(A)$.

Let $A_1=[a_0^{\ell-1}]$ and $A_2$ be two cylinder sets of same length $\ell$. Applying the first part of the lemma, we obtain
$$
\mu(A_1) = \sum_{a_{-p}^{-1},a_{\ell}^{\ell+p-1} :\ a_{-p}^{\ell+p-1}\in \mathcal{L}_{\ell+2p}} \mu([a_{-p}^{\ell+p-1}]) \le |\AA|^{2p} \mu(A_2),
$$
which proves (ii).

It follows that the entropy of $\xi$ is
$$ h(\xi) = \lim_{\ell\to\infty} \dfrac{1}{\ell} \log |\mathcal{L}_\ell|, $$
which clearly achieves the maximum possible entropy for invariant measures with the same support as $\mu$.
\end{proof}

\begin{proof}[Proof of Theorem~\ref{positive entropy}]

Let us assume that $\xi$ is a zero-entropy process satisfying the PIJ property, with the 3-fold self-joining $\lambda$ of width $p$. We have to prove that $\xi$ is a periodic process. 
Without loss of generality, we can assume $p>0$. 

We know that $\xi$ is quasi-uniform. Let $K>0$ be a constant satisfying \eqref{def_quasi_uniform} for all cylinders $B$ and $C$ of the same length. Observe that, if $X$ is a random variable taking values in $\AA$ such that 
\begin{multline}
\label{Q-U for X}
\forall a,b\in\AA, \\\PP(X=a)\PP(X=b)>0\Longrightarrow
\dfrac{1}{K}\PP(X=a)\le \PP(X=b)\le K\PP(X=a)
\end{multline}
then there exists $\delta>0$ satisfying 
\begin{equation}
\label{delta}
H(X)<\delta \Longrightarrow \exists a\in\AA,\ \PP(X=a)=1.
\end{equation}
We now consider the distribution of $\xi_0$ conditioned on $\xi_{-1}, \dots, \xi_{-k}$, $k\ge 1$. It obviously satisfies \eqref{Q-U for X} for any choice of $\xi_{-1}, \dots, \xi_{-k}$. Moreover, since $\xi$ has zero entropy, the measure of the set of sequences $(\xi_{-1}, \dots, \xi_{-k})$ such that 
$H(\xi_0|\xi_{-1}, \dots, \xi_{-k}) < \delta$ tends to 1 as $k\to\infty$:
$$
\mu\left( H(\xi_0|\xi_{-1}, \dots, \xi_{-k}) < \delta \right) \tend{k}{\infty} 1\, .
$$
Therefore
$$
\mu\left( \left\{  [a_{-k}^{-1}] : \xi_0 \mbox{ is determined by } (\xi_{-k}^{-1}= a_{-k}^{-1})\right\} \right) \tend{k}{\infty} 1\, ,
$$
where we say that a random variable $X$ is \emph{determined} by $(Y=y)$ if there exists a value $a$ such that 
$\PP(X=a|Y=y)=1$. 
Let us define the sets of cylinders
$$
\C_k \egdef \left\{ [a_{-k}^{k}] : \xi_{-k-p}^{k+p} \mbox{ is determined by } (\xi_{-k}^{k}= a_{-k}^{k})\right\}\, .
$$
Observe that all cylinders in $\C_k$ have the same measure $\mu_k$ : Indeed, let $A$ and $B$ be two cylinders in $\C_k$ and $A'$ be the cylinder of length $(2k+1)+2p$ such that 
$$
\xi_{-k}^{k}\in A\Longleftrightarrow \xi_{-k-p}^{k+p}\in A'\, .
$$
By Lemma~\ref{lemma_unif}, we have $\mu(A)=\mu(A')\le \mu(B)$. By symmetry, we conclude that $\mu(A)=\mu(B)$. 

Moreover, $\mu(\C_k)\to 1$. 
Hence, for $k$ large enough, there exists a cylinder $[a_{-k}^{k}]\in\C_k$ such that the cylinder $[a_{-k-1}^{k+1}]$ determined by $[a_{-k}^{k}]$ is in $\C_{k+1}$. This implies that $\mu_k=\mu_{k+1}$. 
Let us denote by $\mu_\infty >0$ the ultimate value of the $\mu_k$'s.
If $k$ is large enough so that $\mu(\C_k) > 1-\mu_\infty$, then $\mu(\C_k)=1$. 

We conclude that, for $k$ large enough, $\xi_{-k-p}^{k+p}$ is always determined by $\xi_{-k}^{k}$. 
This means that $\xi$ is periodic. 
\end{proof}

\subsection{An algebraic lemma and its consequence on PIJ processes}

We present in this section an algebraic lemma which can be stated in a very simple context. Let $F$ and $F'$ be two finite sets, and $C$ be a subset of $F\times F'$. We are interested in the situation where we can find two probability measures $m$ on $F$ and $m'$ on $F'$ such that, under the product distribution $m\otimes m'$, $\ind{C}$ is independent of each coordinate. In other words, we demand that there exist two finite families of nonnegative numbers $(m(x))_{x\in F}$ and $(m'(x'))_{x'\in F'}$ such that 

\begin{equation}
\label{finiteIndependence}
\begin{cases}
 \displaystyle \sum_{x\in F} m(x) = \sum_{x'\in F'} m'(x')=1\\
\displaystyle \forall x_0\in F,\ \forall x'_0\in F', \quad
\sum_{x: (x,x_0')\in C} m(x)=\sum_{x': (x_0, x')\in C} m'(x')
\end{cases}
\end{equation}

\begin{figure}[h]
	\begin{center}
	\input{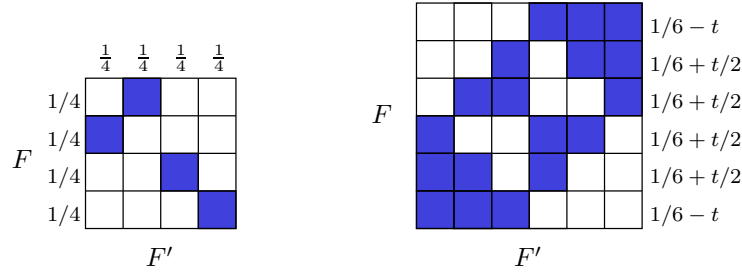}
	\end{center}
\caption{Let $C$ be the colored elements in $F\times F'$. We provide examples of pairs of probability measures satisfying \eqref{finiteIndependence}. On the left, $m\otimes m'(C)=1/4$.
On the right, there is a continuous family of probability measures satisfying \eqref{finiteIndependence} (the measure on $F'$ is defined symmetrically, with a possibly different parameter), but for each pair of them, $m\otimes m'(C)=1/2$. } 
\label{ind_lemma_figure}
\end{figure}

It may happen that many pairs $(m,m')$ satisfy this requirement (see Figure~\ref{ind_lemma_figure}), but as the following lemma shows, $m\otimes m'(C)$ does not depend on the particular choice of the pair $(m,m')$. And the upshot is that this quantity is always a rational number\footnote{It would be interesting to bound the denominator of $m\otimes m'(C)$ by a (possibly linear) function of the sizes of $F$ and $F'$.}!

\begin{lemma}
\label{ind_lemma}
Let $(m_1, m_1')$ and $(m_2, m_2')$ be two pairs of probability measures which satisfy \eqref{finiteIndependence}. Then, $m_1\otimes m_1' (C) = m_2\otimes m_2' (C)\in\QQ$. 
\end{lemma}

\begin{proof}
We fix $x_0\in F$. We have 
$$
m_1\otimes m_1' (C) 
= \sum_{x\in F} m_1(x) \sum_{x':(x,x')\in C} m_1'(x')
= \sum_{x': (x_0, x')\in C} m_1'(x') = m_2\otimes m_1' (C) 
$$
The same argument shows that $m_2\otimes m_1' (C) = m_2\otimes m_2' (C)$. 

Notice that \eqref{finiteIndependence} is a system of linear equations with rational coefficients. Therefore, if it has a solution in $\RR$, it also has a solution in $\QQ$. 
\end{proof}

The preceding lemma has a direct consequence for the process $\xi$ if we assume that it has the PIJ property. Let $\ell\ge1$ be an integer, and denote by $F$ (respectively $F'$) the set of all cylinders $[x_{-p}^{\ell+p-1}]$ of length $\ell+2p$ for $\xi$ (respectively $\xi'$). Then any cylinder $C=[c_{0}^{\ell -1}]$ of length $\ell$ for $\xi''$ can be identified \emph{via} the special pairwise-independent joining $\lambda$ with a subset of $F\times F'$ which is independent of $\xi$ and independent of $\xi'$. Therefore, we get the following corollary.
\begin{corollary}
\label{proba_rationnelle}
If $\xi$ has the PIJ property, then for any cylinder~$C$, 
$$\mu(C)\in\QQ.$$
\end{corollary}

\section{Subshifts of finite type with the PIJ property}
\label{SSFT}
We recall that a \emph{subshift of finite type} is a subset $\Sigma$ of $\AA^\ZZ$ whose elements are all sequences in which a finite number of forbidden patterns never appear. We may always assume that all forbidden patterns have length 2: Replace if necessary $\AA$ by $\AA^{\ell-1}$, where $\ell$ is the largest length of a forbidden word. Then the subshift of finite type can be defined by an $|\AA|$-by-$|\AA|$ matrix $M$, with 0-1 entries, called the \emph{adjacency matrix}: $M_{ab}=1$ if and only if $ab$ is an allowed pattern. 

The subshift of finite type is \emph{irreducible} if $\forall a,b\in\AA$, we can find $k\ge1$ such that $M^k_{ab}\ge1$. In this case, there exists on $\Sigma$ a unique shift-invariant measure of maximal entropy, and the stationary process it defines is quasi-uniform (see~\cite{kitchens1998}, p.~166). This, and the fact that all the examples of processes having the PIJ property that we know arise from subshifts of finite type, justify our interest in the links between the PIJ property and subshifts of finite type.

Observe that there can exist only one quasi-uniform measure whose support is a given irreducible subshift of finite type. Indeed, any such measure would be equivalent to the maximal-entropy measure. Since the latter is ergodic, these measures are equal. Note also that the motivations for the study of the PIJ property allow us to restrict our interest to weakly mixing processes. Therefore we can always assume that the subshift of finite type we consider is aperiodic. 

In the sequel, we only consider subshifts of finite type which are irreducible and aperiodic, endowed with their probability measure of maximal entropy, and use the abbreviation IASFT to denote them. 

\subsection{The Perron value}
Let $\beta$ be the largest eigenvalue of the adjacency matrix $M$, called its \emph{Perron value}. Let $\ell=(\ell_1, \dots, \ell_{|\AA|})$ (respectively $r=(r_1, \dots, r_{|\AA|})$) be a left (respectively right) Perron eigenvector normalized such 
that $\sum \ell_i r_i = 1$. 
We recall (see \cite{kitchens1998} p.~166), that the measure of maximal entropy satisfies, 
for any cylinder $[a_0^k]$, 
\begin{equation}
\label{proba_cylindre}
\mu([a_0^k]) = \frac{\ell_{a_0}r_{a_k}}{\beta^k}\, .
\end{equation}

\begin{theo}
\label{vp_entiere}
Let $\xi$ be an IASFT with the PIJ property. Then the Perron value of the adjacency matrix is an integer. 
\end{theo}
\begin{proof}
Fix $a\in\AA$. We know from Corollary~\ref{proba_rationnelle} that all cylinders have rational probabilities. Hence, $\ell_ar_a\in\QQ$ and $\ell_ar_a/\beta^{k}\in\QQ$ as soon as $M^k_{aa}>0$. 

Since the subshift is aperiodic, $\gcd(\{ k : M^k_{aa} >0\})=1$. 
Therefore, there exist integers $t_1, \dots, t_j$ and $k_1, \dots, k_j$ such that 
$M_{aa}^{k_i}>0$ and $\sum_{i=1}^j t_ik_i = 1$. 
It follows that $\beta^{k_i}\in\QQ$ for all $1\le i\le j$, and 
$\beta^{\sum t_i k_i} = \beta\in\QQ$. 

We let the reader check that if a square matrix with integer coefficients has a rational eigenvalue $\beta$, then $\beta\in\ZZ$. Hence, the Perron value of $M$ is a positive integer. 
\end{proof}

Note that if we drop the aperiodicity assumption, the fact that all cylinders have rational probabilities does not imply that the Perron value is an integer: Consider the subshift of finite type whose adjacency matrix is 
$$
M\egdef \begin{pmatrix}
	0&0&1&0\\
	0&0&0&1\\
	1&1&0&0\\
	1&1&0&0
        \end{pmatrix}\, ,
$$
with Perron value $\sqrt{2}$.

\medskip
Theorem~\ref{vp_entiere} means that, if $\xi$ is an IASFT with the PIJ property, then there exists a positive integer $n$ such that $h(\xi)=\log n$. Since the IASFT we consider is a mixing Markov process, it satisfies the weak-Bernoulli property, and thus is measure-theoretically isomorphic to the full $n$-shift described in Example~\ref{Bernoulli} (see \cite{friedman1970}).

The conclusion of Theorem~\ref{vp_entiere} is also equivalent to the following property (see \cite{lind1995}, p.156): The IASFT is topologically conjugate to another IASFT for which each vertex has $n$ incoming edges and $n$ outcoming edges. We will call such a subshift a \emph{uniform subshift}. 

We may ask whether any IASFT with integer Perron value satisfies the PIJ property. The following section shows that it is not true. From now on, we will denote by $n$ instead of $\beta$ the Perron value of the adjacency matrix $M$.

\subsection{Independence property}

\begin{theo}
\label{independence}
Let $\xi$ be an IASFT with the PIJ property. For $k$ large enough, $\xi_{-\infty}^{0}$ and $\xi_k^{+\infty}$ are independent.
\end{theo}
\begin{proof}
Without loss of generality, we assume that $\xi$ is a uniform subshift. Indeed, the independence property is preserved by topological conjugacy (with a different $k$).
Therefore, for any cylinder $[a_0^k]$, $\mu([a_0^k]) = (|\AA|n^k)^{-1}$. 

We consider a three-fold, pairwise-independent self-joining $\lambda$ of width $p$ such that 
${\xi''}_0^{k} =\varphi({\xi}_{-p}^{k+p},{\xi'}_{-p}^{k+p})$ $\lambda$-a.s. for any $k\ge 0$. 

Let us fix ${\xi}_{-p}^{k+p}$. Then ${\xi''}_0^{k}$ is a function of ${\xi'}_{-p}^{k+p}$. We will first compute the number of preimages of ${\xi''}_0^{k}$. 
\begin{equation}
\label{antecedents}
\forall {\xi}_{-p}^{k+p}, \forall a_0^k, \quad
\left|\left\{ {\xi'}_{-p}^{k+p}: \varphi({\xi}_{-p}^{k+p},{\xi'}_{-p}^{k+p}) = a_0^k \right\}\right| = n^{2p}\, .
\end{equation}
Indeed, for all cylinder $[a_0^k]$, by independence of $\xi$ and $\xi''$, 
$$
\lambda({\xi''}_0^{k} = a_0^k|\xi_{-p}^{k+p}) = \mu({\xi''}_0^{k} = a_0^k) = \frac{1}{|\AA|n^{k}}\, .
$$
On the other hand, this is also equal to 
$$
\sum_{{\xi'}_{-p}^{k+p} : \varphi({\xi}_{-p}^{k+p},{\xi'}_{-p}^{k+p})=a_0^k} \mu([{\xi'}_{-p}^{k+p}])
= \frac{1}{|\AA|n^{k+2p}} \left|\left\{ {\xi'}_{-p}^{k+p}: \varphi({\xi}_{-p}^{k+p},{\xi'}_{-p}^{k+p}) = a_0^k \right\}\right|\, ,
$$
which proves \eqref{antecedents}. 

Recall that $M^k_{a b}$ is the number of cylinders of length $k+1$ starting in $a$ and ending in $b$. Since the IASFT is irreducible aperiodic, there exists $k$ large enough such that $M^{k-2p}_{a b} >0$ for any $a, b\in\AA$. 

Let us fix $a_0$ and $a_k$ in $\AA$. 
Since there are $M^k_{a_0a_k}$ cylinders of length $k+1$ starting in $a_0$ and ending in $a_k$, we deduce from \eqref{antecedents} that, once ${\xi}_{-p}^{k+p}$ is fixed, there are exactly $M^k_{a_0a_k}n^{2p}$ cylinders $[{\xi'}_{-p}^{k+p}]$ yielding ${\xi''}_0 = a_0$ and ${\xi''}_k = a_k$. 
But these cylinders are characterized by $\varphi({\xi}_{-p}^{p},{\xi'}_{-p}^{p}) = a_0$ and 
$\varphi({\xi}_{k-p}^{k+p},{\xi'}_{k-p}^{k+p}) = a_k$. 
Once ${\xi'}_{-p}^{p}$ and ${\xi'}_{k-p}^{k+p}$ are chosen, there are $M^{k-2p}_{{\xi'}_{p}{\xi'}_{k-p}}$ ways of completing the cylinders. 
We thus obtain 
$$
M^k_{a_0a_k}n^{2p} = \sum_{{\xi'}_{-p}^{p}} \sum_{{\xi'}_{k-p}^{k+p}}
M^{k-2p}_{{\xi'}_{p}{\xi'}_{k-p}}
\ind{\{\varphi({\xi}_{-p}^{p},{\xi'}_{-p}^{p})=a_0\}}
\ind{\{\varphi({\xi}_{k-p}^{k+p},{\xi'}_{k-p}^{k+p})=a_k\}} \, .
$$
Observe there are $|\AA|^2n^{4p}$ ways of choosing ${\xi}_{-p}^{p}$ and ${\xi}_{k-p}^{k+p}$, because $M^{k-2p}_{a b} >0$ for any $a, b\in\AA$. 
Therefore, summing over all these possible choices, we get that 
$M^k_{a_0a_k}|\AA|^2n^{6p}$ is equal to 
$$
\sum_{{\xi'}_{-p}^{p}} \sum_{{\xi'}_{k-p}^{k+p}} M^{k-2p}_{{\xi'}_{p}{\xi'}_{k-p}}
\sum_{{\xi}_{-p}^{p}}     \ind{\{\varphi({\xi}_{-p}^{p},{\xi'}_{-p}^{p})=a_0\}}
\sum_{{\xi}_{k-p}^{k+p}}  \ind{\{\varphi({\xi}_{k-p}^{k+p},{\xi'}_{k-p}^{k+p})=a_k\}} \, .
$$
Exchanging the roles of $\xi$ and $\xi'$ in \eqref{antecedents}, we have 
$$
\sum_{{\xi}_{-p}^{p}} \ind{\{\varphi({\xi}_{-p}^{p},{\xi'}_{-p}^{p})=a_0\}} 
= \sum_{{\xi}_{k-p}^{k+p}}  \ind{\{\varphi({\xi}_{k-p}^{k+p},{\xi'}_{k-p}^{k+p})=a_k\}}
= n^{2p}\, .
$$
We conclude that 
$$
M^k_{a_0a_k} = \frac{1}{|\AA|^2n^{2p}} \sum_{{\xi'}_{-p}^{p}} \sum_{{\xi'}_{k-p}^{k+p}} M^{k-2p}_{{\xi'}_{p}{\xi'}_{k-p}}\, ,
$$
which does not depend on $a_0$ and $a_k$. Let us denote by $C_k$ the common value of all coefficients in $M^k$. We can compute $C_k$ by recalling that all cylinders of length $k+1$ have the same mass $(|\AA|n^{k})^{-1}$. Since there are exactly $C_k|\AA|^2$ such cylinders, 
$C_k=n^{k}/|\AA|$. 

Let us see why this implies the independence property for $\xi$. 
For any $j\ge 0$, 
\begin{eqnarray*}
\mu(\xi_{-j}^0=a_{-j}^0, \xi_{k}^{k+j}=a_{k}^{k+j})
&=& C_k\frac{1}{|\AA|n^{k+2j}} \ =\  \left(\frac{1}{|\AA|n^{j}}\right)^2 \\
&=& \mu(\xi_{-j}^0=a_{-j}^0) \mu(\xi_{k}^{k+j}=a_{k}^{k+j})\, .
\end{eqnarray*}
\end{proof}

Remark that since $C_k=n^{k}/|\AA|$ is an integer, we get the following corollary. 
\begin{corollary}
\label{divise}
Let $\xi$ be a uniform IASFT with the PIJ property. Then for $k$ large enough, $|\AA|$ divides $n^k$, where $n$ is the Perron value of $\xi$. 
\end{corollary}

\begin{example}[The triangle subshift, see Figure~\ref{triangle}]
Consider the IASFT on the alphabet $\AA=\{0, 1, 2\}$ with adjacency matrix 
$$
M\egdef \begin{pmatrix}
	0&1&1\\
	1&0&1\\
	1&1&0
        \end{pmatrix}\, .
$$
This is a uniform IASFT whose Perron value is $n=2$. In view of Corollary~\ref{divise}, it does not satisfy the PIJ property.
\begin{figure}[h]
	\begin{center}
	\includegraphics[width=2cm]{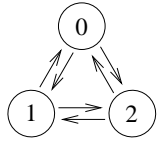}
	\end{center}
\caption{Graph defining the triangle subshift of finite type.}
\label{triangle}
\end{figure}
\end{example}

\subsection{Shift-equivalence to a full shift}
We now wish to investigate whether any IASFT with the PIJ property is topologically conjugate to the full shift on $n$ symbols (where $n$ is the Perron value of the adjacency matrix). Deciding 
whether two IASFT are topologically conjugate is a difficult problem, for which the work of R.F. Williams \cite{williams1973} gives partial answers. Williams introduced purely algebraic properties for square matrices $A$ and $B$ with entries in $\ZZ^+$. 
The matrices $A$ and $B$ are said \emph{strong shift equivalent} if there exist nonnegative integer matrices $R_1,\ldots,R_m$ and $S_1,\ldots,S_m$ with
$$
A = R_1S_1,\quad S_1R_1=R_2S_2,\quad \ldots,\quad S_mR_m=B.
$$
Two IASFT are topologically conjugate if and only if their adjacency matrices are strong shift equivalent. However, strong shift equivalence is not easy to check. That is why the weaker notion of shift equivalence has been introduced:
$A$ is \emph{shift equivalent} to $B$ if there are two nonnegative integer matrices $R$ and $S$ and an integer $m$ such that $RA=BR$, $SB=AS$, $SR=A^m$ and $RS=B^m$. 

\begin{theo}
Let $\xi$ be an IASFT with the PIJ property. Then it is shift-equivalent to the full shift on $n$ symbols, where $n$ is the Perron value of $\xi$.
\end{theo}
\begin{proof}
Without loss of generality, we may still assume that $\xi$ is a uniform subshift. Let $M$ be the adjacency matrix of $\xi$.
We have shown in the proof of Theorem~\ref{independence} that there exists a large enough integer $k$ such that all the entries of $M^k$ are equal. Hence the eigenvalues of $M^k$ are 0 (of order $|\AA|-1$) and $n^k$ (of order~1), and the characteristic polynomial of $M$ is $X^{|\AA|-1}(X-n)$. Lemma~2.2.6 in \cite{kitchens1998}, p.~49 then ensures that $M$ is shift equivalent to the 1-by-1 matrix $(n)$, which itself is shift-equivalent to the $n$-by-$n$ matrix whose entries are all equal to 1.
\end{proof}

Williams had conjectured in 1973 that shift equivalence was the same as topological conjugacy. Kim and Roush have found a counterexample in 1997 \cite{kim1997} which does not involve a full shift. If Williams' conjecture is true when one of the IASFT is the full shift, then any IASFT with the PIJ property is topologically conjugate to the full shift.

\begin{figure}[h]
	\begin{center}
	\includegraphics[width=5cm]{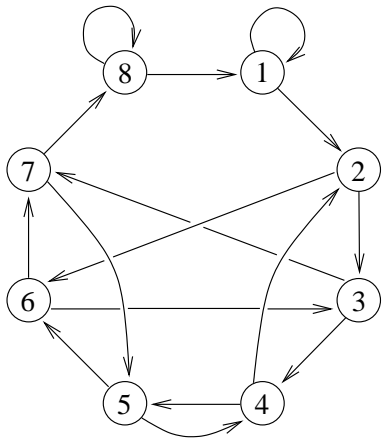}
	\end{center}
\caption{The Ashley subshift.}
\label{ashley}
\end{figure}

Figure~\ref{ashley} shows a graph defining a uniform subshift of finite type which is shift-equivalent to the full 2-shift. This example was proposed by Ashley in 1989. It is not known whether it is topologically conjugate to the full 2-shift (see \cite{kitchens1998}, p.~50). Neither do we know whether this IASFT has the PIJ property.

\section{Does PIJ imply subshift of finite type?}
\label{PIJimplySFT?}

We do not know yet if the PIJ property implies that the process is a subshift of finite type. However, it is true for a slight strengthening of the PIJ property.

\begin{definition}
We say that the process $\xi$ has the \emph{PIJ* property} if there exist a pairwise-independent joining $\lambda$ of $\xi$ and two continuous functions $\varphi_1, \varphi_2:\AA^\ZZ\times\AA^\ZZ\to \AA$ satisfying
\begin{equation}
\label{defPIJeq2}
\xi''_0 = \varphi_1(\xi,\xi') \mbox{ and }\xi_0 = \varphi_2(\xi',\xi'')\qquad(\lambda\mbox{-a.e.})
\end{equation}
\end{definition}

It is easy to find a pairwise-independent joining $\lambda$ satisfying \eqref{defPIJeq} but not \eqref{defPIJeq2}: Take $\xi$ to be the uniform Bernoulli process on $(\ZZ/2\ZZ)^{\ZZ}$, $\xi'$ an independent copy of $\xi$, and set
$$ \xi''_0 = \xi_0+\xi_1+\xi'_0+\xi'_1. $$
However, for this process, we can easily construct another pairwise-independent joining satisfying \eqref{defPIJeq2}. We do not know whether PIJ* is strictly stronger than PIJ.

\begin{theo}
\label{theo:SFT}
Any process $\xi$ enjoying the \emph{PIJ* property} is a subshift of finite type, equipped with an invariant measure of maximal entropy.
\end{theo}

\begin{proof}
We denote by ${\mathcal L}=\cup_\ell{\mathcal L}_\ell$ the set of authorized words in the process $\xi$, that is to say all finite words that occur with positive probability. Observe that $\xi$ is a subshift of finite type if and only if there exists $\ell\in\NN$ such that for any word $W$ of length larger than $\ell$, any $a, b\in\AA$, 
$$aW\in{\mathcal L} \mbox{ and } Wb\in{\mathcal L}\Longrightarrow aWb\in{\mathcal L}.$$

By hypothesis, there exist a pairwise independent joining $\lambda$ of $\xi$ and an integer $p$ such that $\lambda$-almost everywhere, 
\begin{eqnarray}
\xi''_0 &=& \varphi_1(\xi_{-p}^p,{\xi'}_{-p}^p)\label{phi1}\\
\xi_0   &=& \varphi_2({\xi'}_{-p}^p,{\xi''}_{-p}^p)\label{phi2}
\end{eqnarray}

For $U\in{\mathcal L}$ and $d\ge 1$, we define 
$$
\mbox{Ext}(d, U) = \left\{U'\in{\mathcal L}_d: UU'\in{\mathcal L} \right\}.
$$

Observe there exists a word $U\in{\mathcal L}$ of length $|U|>p$ such that for any word $U_1$ with $U_1U\in{\mathcal L}$, $\mbox{Ext}(2p+1, U) = \mbox{Ext}(2p+1, U_1U)$. 
Indeed, it is sufficient to take $U$ such that the cardinal of $\mbox{Ext}(2p+1, U)$ be minimum. 

We consider a word $W$ of length $|W|\ge |U|+2p$, and $a, b\in\AA$ such that $aW\in{\mathcal L}$ and $Wb\in{\mathcal L}$. We have to prove that $aWb\in{\mathcal L}$.

Since $\xi$ and $\xi''$ are independent under $\lambda$, we can observe with positive probability the following situation: 
$\xi_{-|U|+1}^0 = U$ and ${\xi''}_{p-|W|}^p=aW$. 
Therefore, there exist $U_1$ of length $|U_1|=|W|+1-|U|$ and $V$ of length $|V|=|W|+2+2p$ such that with positive probability, we can observe 
$$
\xi_{-|W|}^0 = U_1U , \qquad
{\xi'}_{-|W|}^{2p+1}=V, \qquad
{\xi''}_{p-|W|}^p=aW
$$
all together (see Figure~\ref{fig:sit1} (1)).

\begin{figure}\label{fig:sit1}
\begin{center}
\input{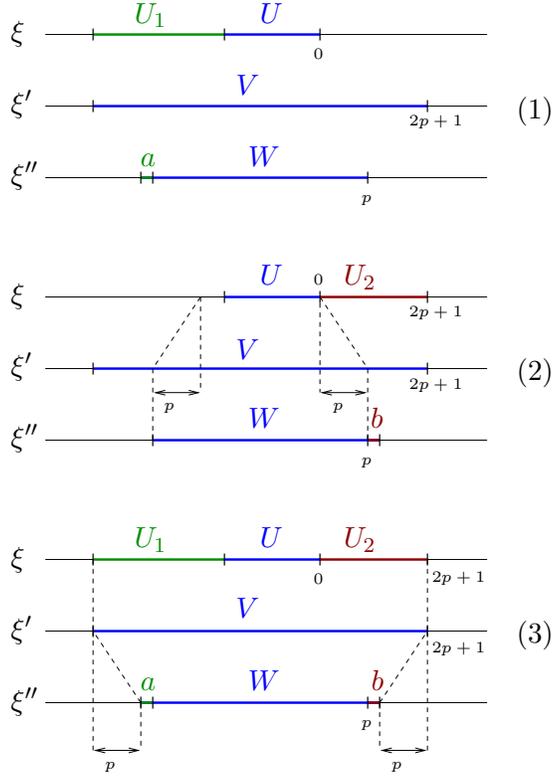}
\end{center}
\caption{The three situations used in the proof of Theorem~\ref{theo:SFT}. In (1), we put together $U$ and $aW$. In (2), we put together $V$ and $Wb$, which forces $U$. In (3), we put together $U_1UU_2$ and $V$, which forces $aWb$.} 
\end{figure}

On the other hand, since $\xi'$ and $\xi''$ are independent under $\lambda$, we can observe with positive probability 
${\xi'}_{-|W|}^{2p+1}=V$ and ${\xi''}_{p-|W|+1}^{p+1}=Wb$. By \eqref{phi2}, this forces $\xi_{-|U|+1}^0=U$. Hence, there exists $U_2\in\mbox{Ext}(2p+1, U)$ such that, with positive probability, we can observe (see Figure~\ref{fig:sit1} (2))
$$
\xi_{-|U|+1}^{2p+1}=UU_2, \qquad
{\xi'}_{-|W|}^{2p+1}=V, \qquad
{\xi''}_{p-|W|+1}^{p+1}=Wb
$$
all together.
But $\mbox{Ext}(2p+1, U) = \mbox{Ext}(2p+1, U_1U)$. Therefore $U_1UU_2\in{\mathcal L}$. 
Since $\xi$ and $\xi'$ are independent under $\lambda$, we can observe with positive probability the following situation (see Figure~\ref{fig:sit1} (3)): 
$\xi_{-|W|}^{2p+1}=U_1UU_2$ and ${\xi'}_{-|W|}^{2p+1}=V$. By \eqref{phi1}, this forces ${\xi''}_{p-|W|}^{p+1}=aWb$. Hence, $aWb\in{\mathcal L}$.

Now, using (iii) of Lemma~\ref{lemma_unif}, we conclude the proof.
\end{proof}

Observe that possessing the PIJ* property does not imply the subshift of finite type is irreducible: Indeed, let $\AA=(\ZZ/2\ZZ)\times(\ZZ/2\ZZ)$, and define the subshift $\Sigma$ as the set of sequences in $\AA^\ZZ$ for which the first coordinate is constant. Let $\mu_0$ and $\mu_1$ be the two ergodic invariant measures on $\Sigma$. This is not an irreducible subshift of finite type, and any convex combination of $\mu_0$ and $\mu_1$ has maximal entropy. Let $\mu=(\mu_0+\mu_1)/2$ be the measure giving the same mass to all cylinders of the same length in $\Sigma$. The process defined by $\mu$ enjoys the PIJ* property with $\xi''_0 = \xi_0+{\xi'}_0$ and $\xi_0  = {\xi'}_0+{\xi''}_0$. 

\section{Further questions}

Our results lead to the conjecture that, if $\xi$ has the PIJ property, then $\xi$ is a continuous coding of one of the two examples presented in Section~\ref{Sec:intro}. If we weaken the requirements of the PIJ property, by replacing the continuity of $\varphi$ with the assumption that $\xi''$ be a \emph{measurable} function of $\xi$ and $\xi'$,  we may ask whether $\xi$ is a measurable coding of one of these examples.

\bibliography{pij.bib}

\providecommand{\bysame}{\leavevmode\hbox to3em{\hrulefill}\thinspace}
\providecommand{\MR}{\relax\ifhmode\unskip\space\fi MR }
\providecommand{\MRhref}[2]{%
  \href{http://www.ams.org/mathscinet-getitem?mr=#1}{#2}
}
\providecommand{\href}[2]{#2}
\begin{thebibliography}{10}

\bibitem{friedman1970}
N.~A. Friedman and D.~S. Ornstein, \emph{On isomorphism of weak {B}ernoulli
  transformations}, Advances in Math. \textbf{5} (1970), 365--394 (1970).

\bibitem{host1991}
B.~Host, \emph{Mixing of all orders and pairwise independent joinings of
  systems with singular spectrum}, Israel J. Math. \textbf{76} (1991), no.~3,
  289--298.

\bibitem{kalikow1984}
S.~A. Kalikow, \emph{Twofold mixing implies threefold mixing for rank one
  transformations}, Ergodic Theory Dynam. Systems \textbf{4} (1984), no.~2,
  237--259.

\bibitem{kim1997}
K.~H. Kim and F.~W. Roush, \emph{The {W}illiams conjecture is false for
  irreducible subshifts}, Electron. Res. Announc. Amer. Math. Soc. \textbf{3}
  (1997), 105--109 (electronic).

\bibitem{kitchens1998}
B.~P. Kitchens, \emph{Symbolic dynamics -- one-sided, two-sided and countable
  state markov shifts}, Universitext, Springer-Verlag, Berlin, 1998.

\bibitem{lind1995}
D.~Lind and B.~Marcus, \emph{An introduction to symbolic dynamics and coding},
  Cambridge University Press, Cambridge, 1995.

\bibitem{delarue2006b}
T.~de~la Rue, \emph{2-fold and 3-fold mixing: why 3-dot-type counterexamples
  are impossible in one dimension}, Bull. Braz. Math. Soc. (N.S.) \textbf{37}
  (2006), no.~4, 503--521.

\bibitem{delarue2006}
\bysame, \emph{An introduction to joinings in ergodic theory}, Discrete Contin.
  Dyn. Syst. \textbf{15} (2006), no.~1, 121--142.

\bibitem{ryzhikov1992}
V.~V. Ryzhikov, \emph{Stochastic wreath products and joinings of dynamical
  systems}, Mat. Zametki \textbf{52} (1992), no.~3, 130--140, 160.

\bibitem{ryzhikov1993}
\bysame, \emph{Joinings and multiple mixing of the actions of finite rank},
  Funktsional. Anal. i Prilozhen. \textbf{27} (1993), no.~2, 63--78, 96.

\bibitem{ryzhikov1993a}
\bysame, \emph{Joinings, wreath products, factors and mixing properties of
  dynamical systems}, Izv. Ross. Akad. Nauk Ser. Mat. \textbf{57} (1993),
  no.~1, 102--128, translation in Russian Acad. Sci. Izv. Math. 42 (1994), no.
  1, 91--114.

\bibitem{williams1973}
R.~F. Williams, \emph{Classification of subshifts of finite type}, Ann. of
  Math. (2) \textbf{98} (1973), 120--153; errata, ibid. (2) 99 (1974),
  380--381.

\end{thebibliography}
\end{document}